\documentclass{article}
\usepackage{amssymb}

\usepackage{graphicx}
\usepackage{amsmath}

\newtheorem{theorem}{Theorem}

\newtheorem{corollary}[theorem]{Corollary}

\newtheorem{definition}[theorem]{Definition}

\newtheorem{proposition}[theorem]{Proposition}

\newenvironment{proof}[1][]{\textbf{Proof.} }{}
\input{tcilatex}

\begin{document}

\title{Nonlinear Connections and Exact Solutions\\
in Einstein and Extra Dimension Gravity}
\author{Sergiu I. Vacaru \\
Departamento de Fisica, Instituto Superior Tecnico, \\
Av. Rovisco Pais 1, Lisbon 1049-001, Portugal,\\
and\\
Departamento de Fisica Aplicada, \\
Facultad de Ciencias Experimentales,\\
Universidad de Huelva, 21071 Huelva, Spain }
\maketitle

\begin{abstract}
We outline a new geometric method of constructing exact solutions of
gravitational field equations parametrized by generic off--diagonal metrics,
anholonomic frames and possessing, in general, nontrivial torsion and
nonmetricity. The formalism of nonlinear connections is elaborated for
(pseudo) Riemannian and Einstein--Cartan--Weyl spaces.
\end{abstract}

\section{ Introduction:\ Nonlinear Connections}

We consider a $(n+m)$--dimensional manifold $V^{n+m},$ provided with general
metric and linear connection structure and of necessary smooth class. It is
supposed that in any point $u\in V^{n+m}$ there is a local splitting into $%
n- $ and $m$--dimensional subspaces, $V_{u}^{n+m}=V_{u}^{n}\oplus
V_{u}^{mn}. $ The local/abstract coordinates are denoted $u=(x,y),$ or $%
u^{\alpha }=\left( x^{i},y^{a}\right) ,$ where $i,j,k,...=1,2,...,n$ and $%
a,b,c,...=n+1,n+2,...,n+m.$ The metric is parametrized in the form 
\begin{equation}
\mathbf{g}=g_{ij}(u)\mathbf{e}^{i}\otimes \mathbf{e}^{j}+h_{ab}(u)\mathbf{e}%
^{a}\otimes \mathbf{e}^{b}  \label{metr}
\end{equation}%
where 
\begin{equation}
\mathbf{\vartheta }^{\mu }=[\vartheta ^{i}=dx^{i},\vartheta
^{a}=dy^{a}+N_{i}^{a}(u)dx^{i}]  \label{ddif}
\end{equation}%
is the dual frame to 
\begin{equation}
\mathbf{e}_{\nu }=[e_{i}=\frac{\partial }{\partial x^{i}}-N_{i}^{a}(u)\frac{%
\partial }{\partial y^{a}},e_{a}=\frac{\partial }{\partial y^{a}}].
\label{dder}
\end{equation}

Let us denote by $\pi ^{\top }:TV^{n+m}\rightarrow TV^{n}$ the differential
of a map $\pi :V^{n+m}\rightarrow V^{n}$ defined by fiber preserving
morphisms of the tangent bundles $TV^{n+m}$ and $TV^{n}.$ The kernel of $\pi
^{\top }$is just the vertical subspace $vV^{n+m}$ with a related inclusion
mapping $i:vV^{n+m}\rightarrow TV^{n+m}.$

\begin{definition}
A nonlinear connection (N--connection) $\mathbf{N}$ on space $V^{n+m}$ is
defined by the splitting on the left of an exact sequence 
\begin{equation*}
0\rightarrow vV^{n+m}\rightarrow TV^{n+m}\rightarrow
TV^{n+m}/vV^{n+m}\rightarrow 0,
\end{equation*}%
i. e. by a morphism of submanifolds $\mathbf{N:\ \ }TV^{n+m}\rightarrow
vV^{n+m}$ such that $\mathbf{N\circ i}$ is the unity in $vV^{n+m}.$
\end{definition}

Equivalently, a N--connection is defined by a Whitney sum of horizontal (h)
subspace $\left( hV^{n+m}\right) $ and vertical (v) subspaces, 
\begin{equation}
TV^{n+m}=hV^{n+m}\oplus vV^{n+m}.  \label{whitney}
\end{equation}%
A space provided with N--connection structure will be denoted $V_{N}^{n+m}.$
We shall use boldfaced indices for the geometric objects adapted to
N--connections. The well known class of linear connections consists a
particular subclass with the coefficients being linear on $y^{a},$ i. e. $%
N_{i}^{a}(u)=\Gamma _{bj}^{a}(x)y^{b}.$

To any sets $N_{i}^{a}(u),$ we can associate certain anholonomic frames (\ref%
{ddif}) and (\ref{dder}), with associated N--connection structure,
satisfying the anholonomy relations 
\begin{equation*}
\lbrack \mathbf{\vartheta }_{\alpha },\mathbf{\vartheta }_{\beta }]=\mathbf{%
\vartheta }_{\alpha }\mathbf{\vartheta }_{\beta }-\mathbf{\vartheta }_{\beta
}\mathbf{\vartheta }_{\alpha }=W_{\alpha \beta }^{\gamma }\mathbf{\vartheta }%
_{\gamma }
\end{equation*}%
with (antisymmetric) nontrivial anholonomy coefficients $W_{ia}^{b}=\partial
_{a}N_{i}^{b}$ and $W_{ji}^{a}=\Omega _{ij}^{a},$ where $\Omega
_{ij}^{a}=e_{[i}N_{j]}^{a}$ are the coefficients of the N--connection
curvature.

Essentially, the method to be presented in this work is based on the notion
of N--connection and considers a Whitney-like splitting of the tangent
bundle to a manifold into a horizontal (see discussion and bibliography in
Refs. \cite{v1,v2,v3}). Here we emphasize that the geometrical aspects of
the N--connection formalism has been studied since the first papers of E.
Cartan \cite{car1}\ and A. Kawaguchi \cite{ak1,ak2}\ (who used it in
component form for Finsler geometry), then one should be mentioned the so
called Ehressman connection \cite{eh}) and the work of W. Barthel \cite{wb}
where the global definition of N--connection was given. The monograph \cite%
{mhss} consider the N--connection formalism elaborated and applied to
geometry of generalized Finsler--Lagrange and Cartan--Hamilton spaces. There
is a set of contributions by Spanish authors, see, for instance, \cite%
{ml1,fe,af}.

We considered N--connections for Clifford and spinor structures \cite%
{vsp1,vsp2}, on superbundles and (super) string theory \cite{vstr1} as well
in noncommutative geometry and gravity \cite{vnc}. The idea to apply the
N--connections formalism as a new geometric method of constructins exact
solutions in gravity theories was suggested in Refs. \cite{vex1,vex2} and
developed and applied in a number of works, see for instance, Ref. \cite%
{vt,vs1,dv}). This contribution outlines the author's and co--authors'
results.

\section{\protect\bigskip N--distinguished Torsions and Curvatures}

The geometric constructions can be adapted to the N--connection structure:

\begin{definition}
A distinguished connection (d--connection) $\mathbf{D}=\{\mathbf{\Gamma }%
_{\beta \gamma }^{\alpha }\}$ on $V_{N}^{n+m}$ is a linear connection
conserving under parallelism the Whitney sum (\ref{whitney}).
\end{definition}

Any d--connection $\mathbf{D}$ is represented by irreducible h- \
v--components $\mathbf{\Gamma }_{\beta \gamma }^{\alpha }=\left( L_{jk}^{i},%
\widetilde{L}_{bk}^{a},C_{jc}^{i},\widetilde{C}_{bc}^{a}\right) $ stated
with respect to N--elongated frames (\ref{ddif}) and (\ref{dder}). This
defines a N--adapted splitting into h-- and v--covariant derivatives, $%
\mathbf{D}=D^{[h]}+D^{[v]},$ where $D^{[h]}=(L,\widetilde{L})$ and $%
D^{[v]}=(C,\widetilde{C}).$ A d--tensor (distinguished tensor, for instance,
a d--metric like (\ref{metr})) formalism and d--covariant differential and
integral calculus can be elaborated \cite{v1} for spaces provided with
general N--connection, d--connection and d--metric structure and nontrivial
nonmetricity 
\begin{equation*}
\mathcal{Q}_{\alpha \beta }\circeq -\mathbf{Dg}_{\alpha \beta }.
\end{equation*}%
The simplest way is to use N--adapted differential forms like $\mathbf{%
\Gamma }_{\beta }^{\alpha }=\mathbf{\Gamma }_{\beta \gamma }^{\alpha }%
\mathbf{\vartheta }^{\gamma }$ with the coefficients defined with respect to
(\ref{ddif}) and (\ref{dder}).

\begin{theorem}
The torsion $\mathcal{T}^{\alpha }\doteqdot \mathbf{D\vartheta }^{\alpha }=d%
\mathbf{\vartheta }^{\alpha }+\Gamma _{\beta }^{\alpha }\wedge \mathbf{%
\vartheta }^{\beta }$ of a d--connection has the irreducible h- v--
components (d--torsions),%
\begin{eqnarray}
T_{\ jk}^{i} &=&L_{\ [jk]}^{i},\ T_{\ ja}^{i}=-T_{\ aj}^{i}=C_{\ ja}^{i},\
T_{\ ji}^{a}=\Omega _{\ ji}^{a},\   \notag \\
T_{\ bi}^{a} &=&T_{\ ib}^{a}=\frac{\partial N_{i}^{a}}{\partial y^{b}}-L_{\
bi}^{a},\ T_{\ bc}^{a}=C_{\ [bc]}^{a}.  \label{dtors}
\end{eqnarray}

\begin{proof}
By a straightforward calculation we can verify the formulas.
\end{proof}
\end{theorem}

The Levi--Civita linear connection $\nabla =\{^{\nabla }\mathbf{\Gamma }%
_{\beta \gamma }^{\alpha }\},$ with vanishing torsion and nonmetricity, is
not adapted to the global splitting (\ref{whitney}). One holds:

\begin{proposition}
There is a preferred, canonical d--connection structure,$\ \widehat{\mathbf{%
\Gamma }}\mathbf{,}$ $\ $on $V_{N}^{n+m},$ constructed only from the metric
coefficients $[g_{ij},h_{ab},N_{i}^{a}]$ and satisfying the conditions $%
\widehat{\mathbf{Q}}_{\alpha \beta }=0$ and $\widehat{T}_{\ jk}^{i}=0$ and $%
\widehat{T}_{\ bc}^{a}=0.$

\begin{proof}
By straightforward calculations with respect to the N--adapted bases (\ref%
{ddif}) and (\ref{dder}), we can verify that the connection 
\begin{equation}
\widehat{\mathbf{\Gamma }}_{\beta \gamma }^{\alpha }=\ ^{\nabla }\mathbf{%
\Gamma }_{\beta \gamma }^{\alpha }+\ \widehat{\mathbf{P}}_{\beta \gamma
}^{\alpha }  \label{cdc}
\end{equation}%
with the deformation d--tensor 
\begin{equation*}
\widehat{\mathbf{P}}_{\beta \gamma }^{\alpha }=(P_{jk}^{i}=0,P_{bk}^{a}=%
\frac{\partial N_{k}^{a}}{\partial y^{b}},P_{jc}^{i}=-\frac{1}{2}%
g^{ik}\Omega _{\ kj}^{a}h_{ca},P_{bc}^{a}=0)
\end{equation*}%
satisfies the conditions of this Proposition. It should be noted that, in
general, the components $\widehat{T}_{\ ja}^{i},\ \widehat{T}_{\ ji}^{a}$
and $\widehat{T}_{\ bi}^{a}$ are not zero. This is an anholonomic frame (or,
equivalently, off--diagonal metric) frame effect.
\end{proof}
\end{proposition}

The torsion of the connection (\ref{cdc}) is denoted $\widehat{\mathbf{T}}%
_{\beta \gamma }^{\alpha }.$ In a similar form we can prove:

\begin{theorem}
The curvature $\mathcal{R}_{\ \beta }^{\alpha }\doteqdot \mathbf{D\Gamma }%
_{\beta }^{\alpha }=d\mathbf{\Gamma }_{\beta }^{\alpha }-\mathbf{\Gamma }%
_{\beta }^{\gamma }\wedge \mathbf{\Gamma }_{\gamma }^{\alpha }$ of a
d--connection $\mathbf{\Gamma }_{\gamma }^{\alpha }$ has the irreducible h-
v-- components (d--curvatures),%
\begin{eqnarray}
R_{\ hjk}^{i} &=&e_{k}L_{\ hj}^{i}-e_{j}L_{\ hk}^{i}+L_{\ hj}^{m}L_{\
mk}^{i}-L_{\ hk}^{m}L_{\ mj}^{i}-C_{\ ha}^{i}\Omega _{\ kj}^{a},  \notag \\
R_{\ bjk}^{a} &=&e_{k}L_{\ bj}^{a}-e_{j}L_{\ bk}^{a}+L_{\ bj}^{c}L_{\
ck}^{a}-L_{\ bk}^{c}L_{\ cj}^{a}-C_{\ bc}^{a}\Omega _{\ kj}^{c},  \notag \\
R_{\ jka}^{i} &=&e_{a}L_{\ jk}^{i}-D_{k}C_{\ ja}^{i}+C_{\ jb}^{i}T_{\
ka}^{b},  \label{dcurv} \\
R_{\ bka}^{c} &=&e_{a}L_{\ bk}^{c}-D_{k}C_{\ ba}^{c}+C_{\ bd}^{c}T_{\
ka}^{c},  \notag \\
R_{\ jbc}^{i} &=&e_{c}C_{\ jb}^{i}-e_{b}C_{\ jc}^{i}+C_{\ jb}^{h}C_{\
hc}^{i}-C_{\ jc}^{h}C_{\ hb}^{i},  \notag \\
R_{\ bcd}^{a} &=&e_{d}C_{\ bc}^{a}-e_{c}C_{\ bd}^{a}+C_{\ bc}^{e}C_{\
ed}^{a}-C_{\ bd}^{e}C_{\ ec}^{a}.  \notag
\end{eqnarray}
\end{theorem}

Contracting the components of (\ref{dcurv}) we prove:

\begin{corollary}
a) The Ricci d--tensor $\mathbf{R}_{\alpha \beta }\doteqdot \mathbf{R}_{\
\alpha \beta \tau }^{\tau }$ has the irreducible h- v--components%
\begin{equation}
R_{ij}\doteqdot R_{\ ijk}^{k},\ \ R_{ia}\doteqdot -R_{\ ika}^{k},\
R_{ai}\doteqdot R_{\ aib}^{b},\ R_{ab}\doteqdot R_{\ abc}^{c}.
\label{dricci}
\end{equation}%
b) The scalar curvature of a d--connection is 
\begin{equation*}
\overleftarrow{\mathbf{R}}\doteqdot \mathbf{g}^{\alpha \beta }\mathbf{R}%
_{\alpha \beta }=g^{ij}R_{ij}+h^{ab}R_{ab}.
\end{equation*}%
c) The Einstein d--densor is computed $\mathbf{G}_{\alpha \beta }=\mathbf{R}%
_{\alpha \beta }-\frac{1}{2}\mathbf{g}_{\alpha \beta }\overleftarrow{\mathbf{%
R}}.$
\end{corollary}

In modern gravity theories one considers more general linear connections
generated by deformations of type $\mathbf{\Gamma }_{\beta \gamma }^{\alpha
}=\widehat{\mathbf{\Gamma }}_{\beta \gamma }^{\alpha }+\mathbf{P}_{\beta
\gamma }^{\alpha }.$We can split all geometric objects into canonical and
post-canonical pieces which results in N--adapted geometric constructions.
For instance, 
\begin{equation}
\mathcal{R}_{\ \beta }^{\alpha }=\widehat{\mathcal{R}}_{\ \beta }^{\alpha }+%
\mathbf{D}\mathcal{P}_{\ \beta }^{\alpha }+\mathcal{P}_{\ \gamma }^{\alpha
}\wedge \mathcal{P}_{\ \beta }^{\gamma }  \label{deformcurv}
\end{equation}%
for $\mathcal{P}_{\beta }^{\alpha }=\mathbf{P}_{\beta \gamma }^{\alpha }%
\mathbf{\vartheta }^{\gamma }.$

\section{Anholonomic Frames and Nonmetricity\newline
in String Gravity}

For simplicity, we investigate here a class of spacetimes when the
nonmetricity and torsion have nontrivial components of type 
\begin{equation}
\mathcal{T}\doteqdot \mathbf{e}_{\alpha }\rfloor \mathcal{T}^{\alpha
}=\kappa _{0}\phi ,\ \mathcal{Q}\doteqdot \frac{1}{4}\mathbf{g}^{\alpha
\beta }\mathcal{Q}_{\alpha \beta }=\kappa _{1}\phi ,\ \mathbf{\Lambda }%
\doteqdot \mathbf{\vartheta }^{\alpha }\mathbf{e}^{\beta }\rfloor (\mathcal{Q%
}_{\alpha \beta }-\mathcal{Q}\mathbf{g}_{\alpha \beta })=\kappa _{2}\phi
\label{aux01}
\end{equation}%
where $\kappa _{0},\kappa _{1},\kappa _{2}=const$ and $\phi =\phi _{\alpha }%
\mathbf{\vartheta }^{\alpha }.$ The abstract indices in (\ref{aux01}) are
''upped'' and ''lowed'' by using $\eta _{\alpha 
{\acute{}}%
\beta 
{\acute{}}%
}$ and its inverse defined from the vielbein decompositions of d--metric, $%
\mathbf{g}_{\alpha \beta }=\mathbf{e}_{\alpha }^{\alpha ^{\prime }}\mathbf{e}%
_{\beta }^{\beta ^{\prime }}\eta _{\alpha ^{\prime }\beta ^{\prime }}.$

Let us consider the strengths $\mathbf{H}_{\nu \mu }\doteqdot \widehat{%
\mathbf{D}}_{\nu }\phi _{\mu }-\widehat{\mathbf{D}}_{\mu }\phi _{\nu
}+W_{\mu \nu }^{\gamma }\phi _{\gamma }$ (intensity of $\phi _{\gamma })$
and $\widehat{\mathbf{H}}_{\nu \lambda \rho }\doteqdot \mathbf{e}_{\nu }%
\mathbf{B}_{\lambda \rho }+\mathbf{e}_{\rho }\mathbf{B}_{\nu \lambda }+%
\mathbf{e}_{\lambda }\mathbf{B}_{\rho \nu }$ (antysimmetric torsion of the $%
\mathbf{B}_{\rho \nu }=-\mathbf{B}_{\nu \rho }$ from the bosonic model of
string theory with dilaton field $\Phi $) and introduce 
\begin{eqnarray*}
\mathbf{H}_{\nu \lambda \rho } &\doteqdot &\widehat{\mathbf{Z}}_{\nu \lambda
\rho }+\widehat{\mathbf{H}}_{\nu \lambda \rho }, \\
\widehat{\mathbf{Z}}_{\nu \lambda } &\doteqdot &\widehat{\mathbf{Z}}_{\nu
\lambda \rho }\mathbf{\vartheta }^{\rho }=\mathbf{e}_{\lambda }\rfloor 
\widehat{\mathbf{T}}_{\nu }-\mathbf{e}_{\nu }\rfloor \widehat{\mathbf{T}}%
_{\lambda }+\frac{1}{2}(\mathbf{e}_{\nu }\rfloor \mathbf{e}_{\lambda
}\rfloor \widehat{\mathbf{T}}_{\lambda })\mathbf{\vartheta }^{\gamma }.
\end{eqnarray*}%
We denote the energy--momentums of fields: 
\begin{equation*}
\mathbf{\Sigma }_{\alpha \beta }^{[\phi ]}\doteqdot \mathbf{H}_{\alpha }^{\
\mu }\mathbf{H}_{\beta \mu }-\frac{1}{4}\mathbf{g}_{\alpha \beta }\mathbf{H}%
^{\nu \mu }\mathbf{H}_{\nu \mu }+\mu ^{2}(\phi _{\alpha }\phi _{\beta }-%
\frac{1}{2}\mathbf{g}_{\alpha \beta }\phi _{\nu }\phi ^{\nu }),
\end{equation*}%
$\mu ^{2}=const,\ \mathbf{\Sigma }_{\alpha \beta }^{[mat]}$ is the source
from any possible matter fields and $\mathbf{\Sigma }_{\alpha \beta }^{[T]}(%
\widehat{\mathbf{T}}_{\nu },\Phi )$ contains contributions of torsion and
dilatonic fields.

\begin{theorem}
\label{theorfields}The dynamics of sigma model of bosonic string gravity
with generic off--diagonal metrics, effective matter and torsion and
nonmetricity (\ref{aux01}) is defined by the system of field equations 
\begin{eqnarray}
\widehat{\mathbf{R}}_{\alpha \beta }-\frac{1}{2}\mathbf{g}_{\alpha \beta }%
\overleftarrow{\mathbf{\hat{R}}} &=&k(\mathbf{\Sigma }_{\alpha \beta
}^{[\phi ]}+\mathbf{\Sigma }_{\alpha \beta }^{[mat]}+\mathbf{\Sigma }%
_{\alpha \beta }^{[T]}),  \label{eecdc1} \\
\widehat{\mathbf{D}}_{\nu }\mathbf{H}^{\nu \mu } &=&\mu ^{2}\phi ^{\mu },\ 
\widehat{\mathbf{D}}^{\nu }(\mathbf{H}_{\nu \lambda \rho })=0,  \notag
\end{eqnarray}%
where $k=const,$ and the Euler--Lagrange equations for the matter fields are
considered on background $V_{N}^{n+m}$).

\begin{proof}
See details in Ref. \cite{v2}.
\end{proof}
\end{theorem}

In terms of differential forms the eqs. (\ref{eecdc1}) are written 
\begin{equation}
\eta _{\alpha \beta \gamma }\wedge \widehat{\mathcal{R}}_{\ }^{\beta \gamma
}=\widehat{\Upsilon }_{\alpha },  \label{eecdc2}
\end{equation}%
where, for the volume 4--form $\eta \doteqdot \ast 1$ with the Hodje
operator ''$\ast $'', $\eta _{\alpha }\doteqdot \mathbf{e}_{\alpha }\rfloor
\eta ,$ $\eta _{\alpha \beta }\doteqdot \mathbf{e}_{\beta }\rfloor \eta
_{\alpha },$ $\eta _{\alpha \beta \gamma }\doteqdot \mathbf{e}_{\gamma
}\rfloor \eta _{\alpha \beta },...,\widehat{\mathcal{R}}_{\ }^{\beta \gamma
} $ is the curvature 2--form and $\Upsilon _{\alpha }$ denote all possible
sources defined by using the canonical d--connection. The deformation of
connection (\ref{cdc}) defines a deformation of the curvature tensor of type
(\ref{deformcurv}) but with respect to the curvature of the Levi--Civita
connection, $\ ^{\nabla }\mathcal{R}_{\ }^{\beta \gamma }.$ The
gravitational field equations (\ref{eecdc2}) transforms into 
\begin{equation}
\eta _{\alpha \beta \gamma }\wedge \ ^{\nabla }\mathcal{R}_{\ }^{\beta
\gamma }+\eta _{\alpha \beta \gamma }\wedge \ ^{\nabla }\mathcal{Z}_{\
}^{\beta \gamma }=\widehat{\Upsilon }_{\alpha },  \label{eecdc3}
\end{equation}%
where $^{\nabla }\mathcal{Z}_{\ \ \gamma }^{\beta }=\nabla \mathcal{P}_{\ \
\gamma }^{\beta }+\mathcal{P}_{\ \ \alpha }^{\beta }\wedge \mathcal{P}_{\ \
\gamma }^{\alpha }.$

\begin{corollary}
A subclass of solutions of the gravitational field equations for the
canonical d--connection defines also solutions of the Einstein equations for
the Levi--Civita connection if and only if $\eta _{\alpha \beta \gamma
}\wedge \ ^{\nabla }\mathcal{Z}_{\ }^{\beta \gamma }=0$ and $\widehat{%
\Upsilon }_{\alpha }=\ ^{\nabla }\Upsilon _{\alpha },$ (i. e. the effective
source is the same for both type of connections).
\end{corollary}

\begin{proof}
It follows from the Theorem \ref{theorfields}.
\end{proof}

This property is very important for constructing exact solutions in Einstein
and string gravity, parametrized by generic off--diagonal metrics and
ahnolonomic frames with associated N--connection structure (see Refs. in %
\cite{v1,v2} and \cite{v3}) and equations (\ref{eecdc3})).

Let us consider a five dimensional ansatz for the metric (\ref{metr}) and
frame (\ref{ddif}) when $u^{\alpha }=(x^{i},y^{4}=v,y^{5});i=1,2,3$ and the
coefficients 
\begin{eqnarray}
g_{ij} &=&diag[g_{1}=\pm
1,g_{2}(x^{2},x^{3}),g_{3}(x^{2},x^{3})],h_{ab}=diag[h_{4}(x^{k},v),h_{5}(x^{k},v)],
\notag \\
N_{i}^{4} &=&w_{i}(x^{k},v),N_{i}^{5}=n_{i}(x^{k},v)  \label{ansatz1}
\end{eqnarray}%
are some functions of necessary smooth class. The partial derivative are
briefly denoted $a^{\bullet }=\partial a/\partial x^{2},a^{^{\prime
}}=\partial a/\partial x^{3},a^{\ast }=\partial a/\partial v.$

\section{Main results:}

\begin{theorem}
\label{ricci}The nontrivial components of the Ricci d--tensors (\ref{dricci}%
) for the ca\-no\-nical d--connection (\ref{cdc}) are 
\begin{eqnarray}
R_{2}^{2} &=&R_{3}^{3}=-\frac{1}{2g_{2}g_{3}}[g_{3}^{\bullet \bullet }-\frac{%
g_{2}^{\bullet }g_{3}^{\bullet }}{2g_{2}}-\frac{(g_{3}^{\bullet })^{2}}{%
2g_{3}}+g_{2}^{^{\prime \prime }}-\frac{g_{2}^{\prime }g_{3}^{\prime }}{%
2g_{3}}-\frac{(g_{2}^{^{\prime }})^{2}}{2g_{2}}],  \notag \\
R_{4}^{4} &=&R_{5}^{5}=-\frac{1}{2h_{4}h_{5}}[h_{5}^{\ast \ast }-h_{5}^{\ast
}(\ln \left| \sqrt{\left| h_{4}h_{5}\right| }\right| )^{\ast }],
\label{riccia} \\
R_{4i} &=&-w_{i}\frac{\beta }{2h_{5}}-\frac{\alpha _{i}}{2h_{5}},\ R_{5i}=-%
\frac{h_{5}}{2h_{4}}[n_{i}^{\ast \ast }+\gamma n_{i}^{\ast }],  \notag
\end{eqnarray}
\begin{equation*}
\alpha _{i}=\partial _{i}h_{5}^{\ast }-h_{5}^{\ast }\partial _{i}\ln \left| 
\sqrt{\left| h_{4}h_{5}\right| }\right| ,\ \beta =h_{5}^{\ast \ast
}-h_{5}^{\ast }[\ln \left| \sqrt{\left| h_{4}h_{5}\right| }\right| ]^{\ast
},\ \gamma =3h_{5}^{\ast }/2h_{5}-h_{4}^{\ast }/h_{4}
\end{equation*}%
$h_{4}^{\ast }\neq 0$ and $h_{5}^{\ast }\neq 0.$
\end{theorem}

\begin{proof}
It is provided in Ref. \cite{v2}.
\end{proof}

\begin{corollary}
The Einstein equations (\ref{eecdc2}) for the ansatz (\ref{ansatz1}) are
compatible for vanishing sources and if and only if the nontricial
components of the source, with respect to the frames (\ref{dder}) and (\ref%
{ddif}), are any functions of type 
\begin{equation*}
\widehat{\Upsilon }_{2}^{2}=\widehat{\Upsilon }_{3}^{3}=\Upsilon
_{2}(x^{2},x^{3},v),\ \widehat{\Upsilon }_{4}^{4}=\widehat{\Upsilon }%
_{5}^{5}=\Upsilon _{4}(x^{2},x^{3})\text{\mbox{ and }}\widehat{\Upsilon }%
_{1}^{1}=\Upsilon _{2}+\Upsilon _{4}.
\end{equation*}

\begin{proof}
The proof, see details in \cite{v2}, follows from the Theorem \ref{ricci}
with the nontrivial components of the Einstein d-tensor, $\widehat{\mathbf{G}%
}_{\ \beta }^{\alpha }=\widehat{\mathbf{R}}_{\ \beta }^{\alpha }-\frac{1}{2}%
\delta _{\ \beta }^{\alpha }\overleftarrow{\mathbf{\hat{R}}},$ computed to
satisfy the conditions 
\begin{equation*}
G_{1}^{1}=-(R_{2}^{2}+R_{4}^{4}),G_{2}^{2}=G_{3}^{3}=-R_{4}^{4}(x^{2},x^{3},v),G_{4}^{4}=G_{5}^{5}=-R_{2}^{2}(x^{2},x^{3}).
\end{equation*}
\end{proof}
\end{corollary}

Having the values (\ref{riccia}), we can proove \cite{v2} the

\begin{theorem}
The system of gravitational field equations (\ref{eecdc1}) (equivalently, of
(\ref{eecdc2})) for the ansatz (\ref{ansatz1}) can be solved in general form
if there are given certain values of functions $g_{2}(x^{2},x^{3})$ (or,
inversely, $g_{3}(x^{2},x^{3})$), $h_{4}(x^{i},v)$ (or, inversely, $%
h_{5}(x^{i},v)$) and of sources $\Upsilon _{2}(x^{2},x^{3},v)$ and $\Upsilon
_{4}(x^{2},x^{3}).$
\end{theorem}

Finally, we note that we have elaborated a new geometric method of
constructing exact solutions in extra dimension gravity and general
relativity theories. The classes of solutions define very general integral
varieties of the vacuum and nonvacuum Einstein equations, in general, with
torsion and nonmetricity and corrections from string theory, and/or
noncommutative/quantum variables. For instance, in five dimensions, the
metrics are generic off--diagonal and depend on four coordinates. So, we
have proved in explicit form how it is possible to solve the Einstein
equations on nonholonomic manifolds (see mathematical problems analized in
Refs. \cite{leit1,leit2}), in our case, provided with nonlinear connection
structure.

{\vskip 4pt}

\textbf{Acknowledgement: }The author is grateful to the Organizes for
supporting his participation at the Conference ''Geometr\'{\i}a de Lorentz,
Murcia 2003, 12-14 de Noviembre'', for partial support from ''Consejeria de
Educacion de la Junta de Andalucia'' and the ''Grupo de Estructura de la
Materia'' of the University of Huelva for his kind hospitality. He also
thanks the referee for valuable remarks and D. Leites, M. de Leon and F.
Etayo for important discussions on the geometry of nonholonomic manifolds
and nonlinear connection formalism.

\end{document}